\documentclass[12pt]{article}
\begin{document}
\newpage
\pagestyle{empty}
\setcounter{page}{0}
%
\newfont{\twelvemsb}{msbm10 scaled\magstep1}
\newfont{\eightmsb}{msbm8}
\newfont{\sixmsb}{msbm6}
\newfam\msbfam
\textfont\msbfam=\twelvemsb
\scriptfont\msbfam=\eightmsb
\scriptscriptfont\msbfam=\sixmsb
\catcode`\@=11
\def\Bbb{\ifmmode\let\next\Bbb@\else
  \def\next{\errmessage{Use \string\Bbb\space only in math mode}}\fi\next}
\def\Bbb@#1{{\Bbb@@{#1}}}
\def\Bbb@@#1{\fam\msbfam#1}
\newfont{\twelvegoth}{eufm10 scaled\magstep1}
\newfont{\tengoth}{eufm10}
\newfont{\eightgoth}{eufm8}
\newfont{\sixgoth}{eufm6}
\newfam\gothfam
\textfont\gothfam=\twelvegoth
\scriptfont\gothfam=\eightgoth
\scriptscriptfont\gothfam=\sixgoth
\def\frak{\frak@}
\def\frak@#1{{\fam\gothfam{{#1}}}}
\def\frak@@#1{\fam\gothfam#1}
\catcode`@=12
%
\def\CC{{\Bbb C}}
\def\NN{{\Bbb N}}
\def\QQ{{\Bbb Q}}
\def\RR{{\Bbb R}}
\def\ZZ{{\Bbb Z}}
\def\cA{{\cal A}}          \def\cB{{\cal B}}          \def\cC{{\cal C}}
\def\cD{{\cal D}}          \def\cE{{\cal E}}          \def\cF{{\cal F}}
\def\cG{{\cal G}}          \def\cH{{\cal H}}          \def\cI{{\cal I}}
\def\cJ{{\cal J}}          \def\cK{{\cal K}}          \def\cL{{\cal L}} 
\def\cM{{\cal M}}          \def\cN{{\cal N}}          \def\cO{{\cal O}}
\def\cP{{\cal P}}          \def\cQ{{\cal Q}}          \def\cR{{\cal R}} 
\def\cS{{\cal S}}          \def\cT{{\cal T}}          \def\cU{{\cal U}}
\def\cV{{\cal V}}          \def\cW{{\cal W}}          \def\cX{{\cal X}}
\def\cY{{\cal Y}}          \def\cZ{{\cal Z}}
\def\refname{On the literature}
\def\qed{\hfill \rule{5pt}{5pt}}
\def\id{\mbox{id}}
\def\Sc{S\!c}
\newtheorem{lemma}{Lemma}
\newtheorem{prop}{Proposition}
\newtheorem{theo}{Theorem}
%
%
\newcommand{\norm}[1]{{\protect\normalsize{#1}}}
\newcommand{\LAP}{{\small E}\norm{N}{\large S}{\Large L}%
  {\large A}\norm{P}{\small P}}
\newcommand{\sLAP}{{\scriptsize E}{\footnotesize{N}}{\small S}%
  {\norm L}{\small A}{\footnotesize{P}}{\scriptsize P}}
\def\logolapin{
  \raisebox{-1.2cm}{\epsfbox
    {/lapphp8/keklapp/ragoucy/paper/enslapp.ps}}}
\def\logolight{{\bf {\large E}{\Large N}{\LARGE S}{\huge L}%
    {\LARGE A}{\Large P}{\large P} }}
\def\logoenslapp{\logolight}
%
%
%
\hbox to \hsize{
\hss
\begin{minipage}{5.2cm}
  \begin{center}
    {\bf Groupe d'Annecy\\ \ \\
      Laboratoire d'Annecy-le-Vieux de Physique des Particules}
  \end{center}
\end{minipage}
\hfill
\logoenslapp
\hfill
\begin{minipage}{4.2cm}
  \begin{center}
    {\bf Groupe de Lyon\\ \ \\
      {\'E}cole Normale Sup{\'e}rieure de Lyon}
  \end{center}
\end{minipage}
\hss}

\vspace {.3cm}
\centerline{\rule{12cm}{.42mm}}

\vfill
\vfill
\begin{center}

  {\LARGE {\bf {\sf On Casimir's Ghost }}} \\[1cm]

\vfill

{\large D. Arnaudon$^{\dagger}$,
        M. Bauer$^{\ddagger}$
        and L. Frappat$^{\dagger}$}

\vfill

{\em $^{\dagger}$ \LAP\ 
\footnote{\ URA 1436 du CNRS, associ{\'e}e {\`a} l'E.N.S. de Lyon et {\`a}
  l'Universit{\'e} de Savoie.
 
\indent
~~~ arnaudon@lapp.in2p3.fr,
bauer@spht.saclay.cea.fr, 
frappat@lapp.in2p3.fr.

\indent   ~Partially supported by European Community Contract
  ERBCHRXCT920069.
}, Chemin de Bellevue BP 110,
74941 Annecy-le-Vieux Cedex, France,

\vspace{3mm}

$^\ddagger$ Service de Physique Th{\'e}orique,
  C.E.A. Saclay, F-91191, Gif-sur-Yvette, France.}

\end{center}

\vfill

\begin{abstract}

We define on the universal enveloping superalgebra of $osp(1|2n)$ a
nonstandard adjoint action, endowing it with a module structure. This
allows, in particular, to construct a bosonic operator which
anticommutes with all the fermionic generators and which appears to be
the square root of a certain Casimir operator. 

\end{abstract}

\vfill
\vfill

\rightline{\LAP-A-587/96}
\rightline{q-alg/9605021}
\rightline{May 96}

\newpage
\pagestyle{plain}

\section{On $osp(1|2n)$\label{sect:introduction}}

It appeared recently \cite{ABosp12q} that the study of representations and 
centre of the quantum universal enveloping superalgebra
$\cU_q(osp(1|2))$ was greatly simplified by the use of a particular
operator called the Scasimir. 
This operator, first written in \cite{Lesniewski} is the
$q$-deformation of a classical operator of $\cU(osp(1|2))$
introduced in \cite{PaisRitt, Pinc}. 
This operator has the following
property: it anticommutes with the fermionic generators and commutes
with the bosonic ones, although it has a bosonic character. Moreover,
the Scasimir can be seen as the square root of the quadratic Casimir
element.

The specificities of the superalgebras $osp(1|2n)$ (compared to
the other simple Lie superalgebras) has led us to think that such a
scheme should also hold for $osp(1|2n)$. The aim of this work was
originally to find the explicit expression of the Scasimir, if it
exists. It appeared that 
the existence of a Scasimir is actually explained by a more general
structure inherent to $osp(1|2n)$, built on a nonstandard adjoint
action. This action endows  $\cU(osp(1|2n))$ with a structure of
$osp(1|2n)$-module, in which a one-dimensional submodule corresponds
to the Scasimir. 

The paper is organized as follows: the end of this section is devoted
to some notations. In section \ref{sect:triangle}, we define the
nonstandard action. We also define a remarkable subspace
$\Omega\subset\cU(osp(1|2n))$. It is proved in section
\ref{sect:sub_rep_tissement} that $\Omega$ is left stable by this
action and its decomposition into simple $osp(1|2n)$-submodules is
accomplished. In section \ref{sect:decomp}, we show that
$\cU(osp(1|2n))$ itself is also a direct sum of finite dimensional
$osp(1|2n)$-modules, which is explicitly given in the case of
$osp(1|2)$. Finally, in section \ref{sect:scasimir}, we compute the
expression of the Scasimir.

\vspace{.5cm}

The basic Lie superalgebra $osp(1|2n)$ is defined, in the Cartan--Weyl
basis, by the generators $\sigma_a$ and $\sigma_{ab}$ for 
$1\leq a,b \leq 2n$ and the relations 
\begin{eqnarray}
  \{\sigma_a , \sigma_b \} &=& \sigma_{ab} \nonumber \\ 
  {}[\sigma_a, \sigma_{bc}] &=& - g_{ab} \sigma_c - g_{ac} \sigma_b \nonumber\\
  {}[\sigma_{ab}, \sigma_{cd}] &=& - g_{ac} \sigma_{bd} - g_{ad} \sigma_{bc}
  - g_{bc} \sigma_{ad} - g_{bd} \sigma_{ac}
  \label{eq:osprel}
\end{eqnarray}
where the $2n\times 2n$ matrix $(g_{ab})$ is given by 
\begin{equation}
  (g_{ab}) = \left(
  \begin{array}{cc}
    0 & -I_n \\
    I_n & 0  \\
  \end{array}
  \right) 
  \label{eq:g_ab}
\end{equation}
where $I_n$ is the $n\times n$ unit matrix.
The bosonic (or even) 
generators $\sigma_{ab}=\sigma_{ba}$ generate the $sp(2n)$
subalgebra of $osp(1|2n)$, while the fermionic (or odd) 
generators $\sigma_a$
form the fundamental representation of $sp(2n)$.
   
A Cartan subalgebra $\cH$ can be obtained by choosing 
$H_a = \sum_b g^{ba}\sigma_{ab}$ 
where $1\leq a \leq n$. 
In the following, we will denote by $\alpha$ the indices in $\{1,\dots,n\}$ 
and define
$\bar\alpha=\alpha + n$ so that $\bar\alpha \in \{n+1,\dots,2n\}$.
The fermionic generators $\sigma_\alpha$ and $\sigma_{\bar\alpha}$
will be called conjugated.
By convention 
$\stackrel{_{\scriptstyle =}}\alpha \;= \alpha$.
The Cartan generators are then the $\sigma_{\alpha\bar\alpha}$. 
With respect to this Cartan subalgebra, the generators $\sigma_\alpha$,
(respectively $\sigma_{\bar\alpha}$),  act as raising (respectively
lowering) operators.
A system of simple root generators can be chosen as 
$\{\sigma_{1\bar 2}, \sigma_{2\bar 3}, \dots, \sigma_{n-1 \bar n},
\sigma_n\}$. Then the positive root generators are 
$\sigma_\alpha$, $\sigma_{\alpha\alpha}$, and
$\sigma_{\alpha\beta}$, $\sigma_{\alpha\bar\beta}$ for
$1\leq\alpha<\beta\leq n$.

\section{On a nonstandard action \label{sect:triangle}}

It is known that the standard adjoint action of the generators
$\sigma_{ab}$ and $\sigma_a$ 
on an element $x$ of $\cU(osp(1|2n))$ given by $[\sigma_{ab},x]$ and
$[\sigma_a,x]_\pm$ 
endows $\cU(osp(1|2n))$ with a structure of
$osp(1|2n)$-module. 

Besides this action, we can define the following alternative action 
\begin{eqnarray}
  \sigma_{ab} \triangleright x &=& [\sigma_{ab}, x] \nonumber\\
  \sigma_a \triangleright x &=& \{\sigma_{a}, x\} \qquad\mbox{for even
    $x$} \nonumber\\
  \sigma_a \triangleright x &=& [\sigma_{a}, x] \,\qquad\mbox{ for odd
    $x$}
  \label{eq:action}
\end{eqnarray}

Obviously, the action of the bosonic part 
$\sigma_{ab}$ through $\triangleright$ again endows 
$\cU(osp(1|2n))$ with a structure of
$sp(2n)$-module, since
\begin{eqnarray}
  {}[\sigma_{ab}, \sigma_{cd}] \triangleright x &=& 
  \sigma_{ab}\triangleright(\sigma_{cd} \triangleright x) 
  - \sigma_{cd}\triangleright(\sigma_{ab} \triangleright x) \;.
  \label{eq:modulebos}
\end{eqnarray}
Moreover, we can check that (\ref{eq:action}) also
defines a structure of $osp(1|2n)$-module because 
$\forall ~ x\in\cU(osp(1|2n))$,
\begin{eqnarray}
  \{\sigma_a, \sigma_b\} \triangleright x &=& 
  \sigma_a\triangleright(\sigma_b \triangleright x) 
  + \sigma_b\triangleright(\sigma_a \triangleright x) \nonumber \\
  {}[\sigma_a, \sigma_{bc}] \triangleright x &=& 
  \sigma_a\triangleright(\sigma_{bc} \triangleright x) 
  - \sigma_{bc}\triangleright(\sigma_a \triangleright x) \;.
  \label{eq:moduleferm}
\end{eqnarray}

\medskip

In the universal enveloping algebra $\cU(osp(1|2n))$, let us define
the following completely antisymmetric polynomials in the fermionic
generators $\sigma_a$:
\begin{equation}
{}[a_1 a_2 \dots a_p] \equiv \sum_{s\in {\frak S}_p} \varepsilon(s) \ 
            \sigma_{a_{s(1)}}  \sigma_{a_{s(2)}} \dots \sigma_{a_{s(p)}}  
  \label{eq:crochet}
\end{equation}
where ${\frak S}_p$ is the permutation group of $p$ elements and
$\varepsilon(s)$ is the signature of the permutation~$s$. Note that
this expression vanishes when two indices $a_i$ and $a_j$ coincide
and that consequently $[a_1 a_2 \dots a_p]=0$ for $p>2n$. 
We will denote by $\Omega$ the vector subspace of  $\cU(osp(1|2n))$
spanned by the $[a_1 a_2 \dots a_p]$ for $p\in\{0,\dots,2n\}$, with
the convention that for $p=0$, the value of $[\dots]$ is $1$.

\begin{prop}
  We have the decomposition
  \begin{equation}
    \cU(osp(1|2n)) = \cU(sp(2n)) \otimes \Omega \;.
    \label{eq:decomp}
  \end{equation}
  \label{propU}
\end{prop}
{\bf Proof:}
A Poincar{\'e}--Birkhoff--Witt basis ${\frak B}$ of $\cU(osp(1|2n))$ 
is given by the tensor product of a Poincar{\'e}--Birkhoff--Witt basis
of $\cU(sp(2n))$ by elements of the form 
$\sigma_{a_1} \dots \sigma_{a_p}$ with $1\leq a_1<\dots<a_p\leq 2n$.  
Considering an element $\omega\in\Omega$, it can be re-expressed in
${\frak B}$ using the commutation relations (\ref{eq:osprel}). The
so-obtained result contains {\em i)} a product of fermions of the same
degree as $\omega$ and {\em ii)} a finite linear combination of elements of
${\frak B}$  with a fermionic part of strictly lower degree. It
follows that there exists an invertible triangular matrix relating the
basis ${\frak B}$ and the tensor product of the 
Poincar{\'e}--Birkhoff--Witt basis of $\cU(sp(2n))$ by the basis of
$\Omega$ given by the elements of the form (\ref{eq:crochet}). 

\qed

\section{On a nice subrepresentation \label{sect:sub_rep_tissement}}

We will now consider the action of $osp(1|2n)$ by $\triangleright$ on
$\Omega$. 

\begin{prop}
  The vector space $\Omega$ is stable under the action of any
  $\sigma_{ab}$ by $\triangleright$, i.e. 
  \begin{eqnarray}
    \sigma_{ab} \triangleright \Omega &\subset& \Omega \;. 
    \label{eq:sigmabOmega}
  \end{eqnarray}
  Then, under this action, $\Omega$ is a $sp(2n)$-submodule of
  $\cU(osp(1|2n))$.  Moreover,  
  $\cU(osp(1|2n))$ and $\cU(sp(2n)) \otimes \Omega$
  are equivalent as  $sp(2n)$-modules. 
  \label{propab}
\end{prop}

\medskip
\noindent
{\bf Proof:}
The commutation relations (\ref{eq:osprel}) and the definition
(\ref{eq:action}) lead to the formula
\begin{equation}
  \sigma_{bc} \triangleright [a_1 a_2 \dots a_p] = 
 - \sum_{i=1}^p  g_{ba_i} [a_1 \dots \stackrel{^{(i)}}{c} \dots a_p]
 - \sum_{i=1}^p  g_{ca_i} [a_1 \dots \stackrel{^{(i)}}{b} \dots a_p]
 \in \Omega
  \label{eq:actionbc}
\end{equation}
in which the notation $\stackrel{^{(i)}}{x}$ means that $a_i$ is
replaced by $x$. 

\qed

\begin{prop}
  The vector space $\Omega$ is stable under the action of any
  $\sigma_a$ by $\triangleright$, i.e. 
  \begin{eqnarray}
    \sigma_a \triangleright \Omega &\subset& \Omega \;.
    \label{eq:sigmaOmega}
  \end{eqnarray}
  Then, the action of $\triangleright$ endows $\Omega$ with a
  structure of $osp(1|2n)$-submodule of $\cU(osp(1|2n))$. 
  \label{propa}
\end{prop}

\medskip\noindent
{\bf Proof:}
The stability statement is proved by a direct calculation. We consider
four cases, depending on whether $a$, $\bar a$ belong to
$\{a_1,\dots,a_p\}$ or not. 

Let us first take the case where $a,\bar a \not\in
\{a_1,\dots,a_p\}$. To this aim, we consider the polynomial $[a a_1
a_2 \dots a_p]$ 
of $\Omega$, in which we move the $\sigma_a$ generator to the left using
the commutation relations (\ref{eq:osprel}). We get
\begin{eqnarray}
  [a a_1 a_2 \dots a_p] &=& \sigma_a \sum_{2k\leq p} ~~ 
  \sum_{a'_1,\dots,a'_{p-2k}} [a'_1 \dots a'_{p-2k}] ~ 2^k k! ~  
  {p+1 \choose 2k+1} \prod_{i=1}^k (-g_{..}) \nonumber\\
  && -  \sum_{j=1}^p \sigma_{a a_j} \sum_{2k\leq p-1} ~~ 
  \sum_{a'_1,\dots,a'_{p-2k-1}} [a'_1 \dots a'_{p-2k-1}] ~ 2^k k! ~  
  {p+1 \choose 2k+2} \prod_{i=1}^k (-g_{..}) \nonumber\\
  \label{eq:case1}
\end{eqnarray}
$k$ counts for the number of pairs of conjugated fermionic generators 
missing in the right hand side, each pair $\sigma_b, \sigma_{\bar b}$ 
being replaced by a factor $-g_{b \bar b}$. The remaining indices $a'_1$,
$a'_2,\dots$ shall appear in the same order as in the left hand side. \\
We repeat the procedure, starting again from $[a a_1 a_2 \dots a_p]$,
by moving the $\sigma_a$ generator to the right. We then symmetrize or
antisymmetrize according to the parity of $p$. One obtains 
\begin{eqnarray}
  2 [a a_1 a_2 \dots a_p] &=& \sum_{2k\leq p} ~~ 
  \sum_{a'_1,\dots,a'_{p-2k}} 
  \sigma_a \triangleright [a'_1 \dots a'_{p-2k}] ~ 2^k k! ~  
  {p+1 \choose 2k+1} \prod_{i=1}^k (-g_{..}) \nonumber\\
  && -  \sum_{j=1}^p ~ \sum_{2k\leq p-1} ~~ 
  \sum_{a'_1,\dots,a'_{p-2k-1}} 
  \sigma_{a a_j} \triangleright [a'_1 \dots a'_{p-2k-1}] ~ 2^k k! ~  
  {p+1 \choose 2k+2} \prod_{i=1}^k (-g_{..}) \nonumber\\
  \label{eq:case1suite}
\end{eqnarray}
Notice that the first term of the right hand side with $k=0$ reduces
to $\sigma_a \triangleright [a_1 a_2 \dots a_p]$. It follows that 
$\sigma_a \triangleright [a_1 a_2 \dots a_p]$ can be expressed in terms
of the left hand side and lower order terms of the right hand side (by
virtue of Proposition \ref{propab}, the second sum is an element of
$\Omega$). Then, by recursion on $p$, $\sigma_a \triangleright [a_1
a_2 \dots a_p]\in \Omega$.

We now consider the case when $a$, but not $\bar a$, belongs to the set
$\{a_1,\dots,a_p\}$. Taking for instance $a=a_1$, we have
\begin{equation}
  \sigma_a \triangleright [a a_2 \dots a_p] = 0
  \label{eq:case2}
\end{equation}

The case when $\bar a$, but not $a$, belongs to the set
$\{a_1,\dots,a_p\}$ can be treated in an analogous way as the first
one. One finds, taking for example $\bar a = a_1$,  
\begin{eqnarray}
  2 [a \bar a a_2 \dots a_p] &=& \sum_{2k\leq p-1} ~~ 
  \sum_{a'_1,\dots,a'_{p-2k}} 
  \sigma_a \triangleright [a'_1 \dots a'_{p-2k}] ~ 2^k k! ~  
  {p+1 \choose 2k+1} \prod_{i=1}^k (-g_{..}) \nonumber\\
  && +  \sum_{j=2}^p ~ \sum_{2k\leq p-2} ~~ 
  \sum_{a'_1,\dots,a'_{p-2k-1}} 
  \sigma_{a a_j} \triangleright [a'_1 \dots a'_{p-2k-1}] ~ 2^k k! ~  
  {p+1 \choose 2k+2} \prod_{i=1}^k (-g_{..}) \nonumber\\
  && + 2 g_{a\bar a} \sum_{2k\leq p-2} ~~ 
  \sum_{a'_1,\dots,a'_{p-2k-1}} 
  [a'_1 \dots a'_{p-2k-1}] ~ 2^k k! ~  
  {p+1 \choose 2k+3} \prod_{i=1}^k (-g_{..}) \;.
  \label{eq:case3}
\end{eqnarray}
As before, this allows us to express 
$\sigma_a \triangleright [\bar a a_2 \dots a_p] $ as a combination of elements
of $\Omega$ and lower order terms. By recursion on $p$, this proves
that $\sigma_a \triangleright [\bar a a_2 \dots a_p]\in \Omega$.

The last case is when $a,\bar a\in \{a_1,\dots,a_p\}$.
Taking $a=a_1$ and $\bar a=a_2$, we now directly compute 
\begin{eqnarray}
  && 2\sigma_a  \triangleright  [a \bar a a_3 \dots a_p] = \nonumber\\
  &&\qquad -g_{a\bar a} \sum_{2k\leq p-2} ~~ 
  \sum_{a'_1,\dots,a'_{p-2k-2}} 
  \sigma_a \triangleright [a'_1 \dots a'_{p-2k-2}] ~ 2^k k! ~ (2k+1) ~   
  {p+1 \choose 2k+3} \prod_{i=1}^k (-g_{..}) \nonumber\\
  &&\qquad + g_{a\bar a}  \sum_{j=3}^p ~ \sum_{2k\leq p-3} ~~ 
  \sum_{a'_1,\dots,a'_{p-2k-3}} 
  \sigma_{a a_j} \triangleright [a'_1 \dots a'_{p-2k-3}] ~ 2^{k+1} (k+1)! ~  
  {p+1 \choose 2k+4} \prod_{i=1}^k (-g_{..})  \nonumber\\
  \label{eq:case4}
\end{eqnarray}

\qed

\begin{prop}
  As an $osp(1|2n)$-module, $\Omega$ is equivalent to the direct sum of
  the representations $\cV_j$ characterized by the highest weights 
  $(\underbrace{1,\dots,1}_j, \underbrace{ 0,\dots,0}_{n-j})$
  with respect to $\sigma_{1\bar 1},\dots,\sigma_{n\bar n}$:
  \begin{equation}
    \Omega \simeq \bigoplus_{j=0}^n ~ \cV_j
    \label{eq:direct_sum}
  \end{equation}
  \label{propdec}
  Notice that $\cV_0$ and $\cV_1$ are the trivial and fundamental
  representations of $osp(1|2n)$. 
\end{prop}

\medskip\noindent
{\bf Proof:} 
Let us recall that all the finite dimensional representations 
of $osp(1|2n)$ are
completely reducible \cite{KacSuper,ScheuSuper}. 
The representation $\Omega$
thus decomposes into a sum of irreducible highest weight 
representations. What we need is to recognize in $\Omega$ the highest
weights of the representations $\cV_j$, which are
irreducible. 
We will show that the highest weight vectors have the form 
\begin{equation}
  v_j = [1 ~ 2 ~ \ldots  ~ j ~ j+1 ~ \overline{j+1} ~ \dots  ~ n ~ \bar n] 
  ~ + ~ \mbox{lower order terms}
  \label{eq:hwv}
\end{equation}
where all the indices but $\bar 1,\dots,\bar j$ are present in the
first term.
The lower order terms are obtained from the first one by
deleting pairs of conjugated indices $\alpha$, $\bar\alpha$. These
highest weight vectors $v_j$ have highest weight 
\begin{equation}
  (\underbrace{1,\dots,1}_j, \underbrace{ 0,\dots,0}_{n-j})
  \label{eq:weights}
\end{equation}
with respect to $\sigma_{1\bar 1},\dots,\sigma_{n\bar n}$. 
We first note that 
$\sigma_{\alpha\;\overline{\alpha+1}} ~ \triangleright v_j = 0$
as long as the lower order terms in $v_j$ are symmetric in the
exchange of pairs of indices $\{\alpha, \bar\alpha\}$ and 
$\{\alpha+1, \overline{\alpha+1}\}$ for $j+1 \leq \alpha \leq n-1$.
We now need to construct $v_j$, with this specification, such that
$\sigma_n\triangleright v_j = 0$. Let us first consider the action of
$\sigma_n$ on the leading term 
$[1 ~ 2 ~ \ldots  ~ j ~ j+1 ~ \overline{j+1} ~ \dots  ~ n ~ \bar n]$. From
(\ref{eq:case4}) and repeated uses of (\ref{eq:actionbc}) and
(\ref{eq:case1suite}), one obtains 
\begin{equation}
  \sigma_n\triangleright [1 ~ 2 ~ \ldots  ~ j ~ j+1 ~ \overline{j+1} ~
  \dots  ~ n ~ \bar n]  = \sigma_n\triangleright \omega
  \label{eq:demo1}
\end{equation}
where $\omega\in\Omega$ is a linear combination of terms obtained from
the leading one by deleting the pair $n$, $\bar n$ and possibly other 
pairs of conjugated indices. By construction, $\omega$ is invariant
under the exchange of pairs of indices $\{\alpha, \bar\alpha\}$ and 
$\{\beta, \bar\beta\}$
for $j+1 \leq \alpha,\beta \leq n-1$. Let $\tilde \omega$ be the
expression obtained from $\omega$ by extending this symmetry to the
pair $\{n, \bar n\}$. One can remark that {\em i)} the degree of
$\omega$ and $\tilde\omega$ is strictly less than the degree of the
leading term and {\em ii)} all the terms of $\tilde\omega -\omega$
contain the pair of indices $\{n,\bar n\}$. It follows that we can
compute $\sigma\triangleright(\tilde\omega-\omega)$ as we computed the
action of $\sigma_n$ on the leading term, and we get
\begin{equation}
  \sigma_n\triangleright \Big([1 ~ 2 ~ \ldots  ~ j ~ j+1 ~
  \overline{j+1} ~ \dots  ~ n ~ \bar n]  - \tilde\omega \Big)=
  \sigma_n\triangleright \omega' 
  \label{eq:demo2}
\end{equation}
where the degree of $\omega'$ is strictly less than the degree of
$\omega$. Then, repeating the whole procedure as many times as
necessary, we finally obtain $\sigma_n \triangleright v_j=0$, where 
the element $v_j$ is 
given by (\ref{eq:hwv}), the lower order terms being completely
symmetric under the exchange of pairs of indices $\{\alpha,
\bar\alpha\}$ and  $\{\beta, \bar\beta\}$ for 
$j+1 \leq \alpha,\beta \leq n$. 
Since $v_j$ is annihilated by all the simple positive root generators
of $osp(1|2n)$, it is a highest weight vector. 

The representation generated by the highest weight vector $v_j$ is
$\cV_j$, of dimension ${2n+1\choose j}$. 
The sum of these dimensions is $2^{2n}=\dim
\Omega$ which concludes the proof. 

\qed

The representations $\cV_j$ involved in this decomposition actually
correspond 
\cite{RittSch} to the $so(2n+1)$ representations obtained as the
external powers $\bigwedge\nolimits^j \cV$ of the fundamental
representation $\cV$ of $so(2n+1)$.

\medskip\noindent
{\bf Corollary:} {\em There exists an element $\Sc$ of
  $\cU(osp(1|2n))$ belonging to  $\Omega$ that
  anticommutes with all the fermionic generators and commutes with all
  the bosonic ones. This element $\Sc$ is even with respect to the
  $\ZZ_2$-gradation of $\cU(osp(1|2n)$ and hence is not a Casimir
  operator. It is rather the square root of a central element. We call
  it the Scasimir. As in \cite{PaisRitt} for $osp(1|2)$, 
  a normalized version of this
  operator would provide a realization of $(-1)^H$, $H=\sum_{a=1}^n
  H_a$ being the principal gradation or the ``fermion number''.
}

The proof of the corollary is obvious taking $\Sc = v_0$, the highest
weight vector of the trivial $osp(1|2n)$-representation of $\Omega$. 
The explicit expression of $\Sc$ is given in section \ref{sect:scasimir}.

\section{On the decomposition of $\cU(osp(1|2n))$ \label{sect:decomp}}

\begin{prop}
  As an $osp(1|2n)$-module (with the module structure endowed by the
  action of $\triangleright$), the universal enveloping superalgebra
  $\cU(osp(1|2n))$ splits into finite dimensional representations. 
  \label{propdec2}
\end{prop}

\medskip\noindent
{\bf Proof:} 

Consider the natural filtration of $\cU(sp(2n))$
\begin{equation}
  \CC = \cU_0 \subset \cU_1 \subset \cU_2 \subset \dots \subset \cU_d
  \subset \dots \subset \cU(sp(2n)) = \bigcup^\infty \cU_d
  \label{eq:filtre}
\end{equation}
allowing us to use the degree $d$ of an element $b \in \cU(sp(2n))$
as the smallest $d$ such that $b \in \cU_d$.
\\
Let us now prove that $\cU_d\otimes\Omega$ is stable under the action of
$\triangleright$, which is enough to achieve the proof of the
proposition. (Notice that this generalizes the proposition \ref{propa}
which states the stability of $\cU_0\otimes\Omega$). 
Taking an element $b\omega \in \cU(osp(1|2n)) = \cU(sp(2n)) \otimes
\Omega$ with $b \in \cU_d$ and $\omega\in\Omega$, one has for any
fermionic generator $\sigma\in osp(1|2n)$ 
\begin{equation}
  \sigma \triangleright b\omega = [\sigma, b] \omega + b \sigma
  \triangleright \omega \;.
  \label{eq:action_sigbom}
\end{equation}
{}From Proposition \ref{propa}, the term $b \sigma \triangleright \omega$
belongs to $\cU_d \otimes \Omega$.
Using the commutation relations, the term $[\sigma, b]$ can be
expressed as a sum of terms of the type $b_a\sigma_a$ with
$b_a\in\cU_{d-1}$. We are then led to check that $\sigma_a \omega \in
\cU_1 \otimes \Omega$. Using 
\begin{equation}
  \sigma_a \omega = \frac12 (\sigma_a \triangleright \omega + [\sigma_a,
  \omega]_{\pm}) 
  \label{eq:prod_action}
\end{equation}
where $[~,~]_{\pm}$ denotes the usual $\ZZ_2$-graded commutator, and
equation (\ref{eq:sigmaOmega}), and writing $\omega$
as in (\ref{eq:crochet}), one is left with
\begin{eqnarray}
  [\sigma_a, \omega]_{\pm} &=& \sum_{s\in {\frak S}_p} \varepsilon(s)
  \sum_{j=1}^p (-1)^{j-1} 
  \sigma_{a_{s(1)}}\dots\sigma_{a_{s(j-1)}} \{\sigma_a,
  \sigma_{a_{s(j)}} \} \sigma_{a_{s(j+1)}}
  \dots \sigma_{a_{s(p)}} \nonumber\\
  &=& \sum_{s\in {\frak S}_p} \varepsilon(s)
  \sum_{j=1}^p (-1)^{j-1} \Big( \sigma_{aa_{s(j)}} \sigma_{a_{s(1)}} \dots
  \sigma_{a_{s(j-1)}}
  \sigma_{a_{s(j+1)}} \dots \sigma_{a_{s(p)}} \nonumber\\
  && + \sigma_{a_{s(1)}} \dots \sigma_{a_{s(k-1)}} 
  (g_{aa_{s(j-1)}} \sigma_{a_{s(j)}} + g_{a_{s(j)} a_{s(j-1)}} \sigma_a) 
  \sigma_{a_{s(k)}} \dots \sigma_{a_{s(j-1)}} \sigma_{a_{s(j+1)}}
  \dots \sigma_{a_{s(p)}} \Big) \nonumber \\[.3cm]
  &&
  \label{eq:com_sigom}
\end{eqnarray}
We use recursively the following procedure to the expression 
(\ref{eq:com_sigom}) : {\em i)} the fermionic
generator $\sigma_a$ is 
put in the first place, using the commutation relations
(\ref{eq:osprel}), {\em ii)} the created bosonic
generators are also put in the first place, {\em iii)} one uses
iteratively equations (\ref{eq:sigmaOmega}), (\ref{eq:prod_action})
and (\ref{eq:com_sigom}) on the generated terms of
lower order. We eventually get a  result that belongs to 
$\cU_1 \otimes \Omega$. 
Therefore $\sigma \triangleright b \omega \in \cU_d \otimes
\Omega$, hence $\cU_d \otimes \Omega$ is stable under the action of
$\triangleright$. 

\qed

We now give the explicit decomposition in the case of $osp(1|2)$. 
We start with the decomposition of $\cU_d\otimes\Omega$ in terms of
$sp(2)$-representations. The action $\triangleright$ and the usual
adjoint action of generators of $sp(2)$ are indeed identical. Since
$\cU_d\otimes\Omega$ is finite dimensional, its decomposition in terms
of $osp(1|2)$-representations is therefore uniquely determined. One
finds 
\begin{equation}
  \cU_d \otimes \Omega = \cR_{d+\frac12} \oplus \cR_d \oplus
  \cR_{d-\frac12} \oplus \dots \oplus \cR_0 \oplus C
  \cU_{d-2}\otimes\Omega 
  \label{eq:decosp}
\end{equation}
where $C$ is the quadratic Casimir of $osp(1|2)$. Hence, under the
action of $\triangleright$, $\cU(osp(1|2))$ decomposes into
$osp(1|2)$-representation as follows:
\begin{equation}
  \cU(osp(1|2)) \simeq \cZ \otimes \bigoplus_{j\in\NN/2} \cR_j
  \label{eq:deccU}
\end{equation}
where $\cZ$ denotes the centre of $\cU(osp(1|2))$, i.e. the algebra
generated by $1$ and $C$. This decomposition is actually the same as
that obtained from the usual adjoint action \cite{Pinc}. 

We conjecture that this result generalizes for the case of
$osp(1|2n)$, i.e.  the decomposition of $\cU(osp(1|2n))$ into
$osp(1|2n)$-representations is the same for the two actions (adjoint
and $\triangleright$).
This is motivated by the following arguments: first we know that the
decompositions in terms of $sp(2n)$-representations coincide (which
would be sufficient for finite dimensional representations). Second
the ambiguity in the gathering of $sp(2n)$-representations into
$osp(1|2n)$-representations, which arises in the infinite dimensional
case, can be lifted by looking at the finite dimensional subspaces
$\cU_d\otimes\Omega$ for small~$d$. 

\section{On the Scasimir operator\label{sect:scasimir}}
\subsection{On its explicit expression}

This section is devoted to the construction of $\Sc=v_0$, the element
of $\Omega$ that anticommutes with all the fermionic generators. As
stated in (\ref{eq:hwv}), we start from 
\begin{equation}
  v_0 = [1 ~ \bar 1 ~ 2 ~ \bar 2 ~ \dots  ~ n ~ \bar n]
  ~ + ~ \mbox{lower order terms.}
  \label{eq:v0}
\end{equation}
We define for convenience
\begin{equation}
  A^{(n)}_{2n-2k} \equiv \sum_{1\leq \alpha_1<\dots < \alpha_k\leq n}
  {}[1 ~ \bar 1 ~ \ldots ~ 
  \widehat{\alpha_1} ~ \widehat{\overline{\alpha_1}} ~ \dots ~ 
  \widehat{\alpha_k} ~ \widehat{\overline{\alpha_k}} ~ \dots ~ 
  n ~ \bar n] 
  \label{eq:A2n}
\end{equation}
in which the hat signals missing indices. 
For instance, $A^{(n)}_0=1$ and
\begin{equation}
  \begin{array}{lll}
    A^{(1)}_2 = [1 \bar 1] &&\\[.3cm]
    A^{(2)}_4 = [1\bar 1 2\bar 2] \quad & A^{(2)}_2 = [1 \bar 1]+[2
    \bar 2] \quad &  \\[.3cm]
    A^{(3)}_6 = [1\bar 1 2\bar 2 3\bar 3] \quad & A^{(3)}_4 = [1 \bar
    1 2\bar 2]+[1\bar 1 3\bar 3]+[2\bar 2 3\bar 3] \quad & A^{(3)}_2 =[1 \bar
    1]+[2\bar 2]+[3\bar 3]  \\
  \end{array}
  \label{eq:A2nexemple}
\end{equation}
The expressions
$A^{(n)}_{2n-2k}$ are invariant with respect to permutations 
of pairs of conjugated indices and they correspond,
according to the previous section, to the terms involved in $v_0$:
\begin{equation}
  v_0 = \sum_{k=0}^n x^{(n)}_k A^{(n)}_{2n-2k}
  \qquad\mbox{with $x^{(n)}_0=1$}\;.
  \label{eq:v0suite}
\end{equation}
{}From the equations used in Proposition \ref{propa}, one can compute
explicitly 
\begin{equation}
  \sigma_a \triangleright v_0 = \sigma_a C_a + \sum_b \sigma_{ab} C_{ab}
  \label{eq:sigmav0}
\end{equation}
Now, demanding that the output is zero, one gets the two following
recursion relations:
\begin{eqnarray}
  &&\!\!\!\!\!\!\!\!\! 
  y_k^{(n)} = -\frac 1 k \sum_{p=0}^{k-1} y_p^{(n)} 
  \frac{(2k+2)!/4}{(2p+2)!(2k-2p+1)!} \Big((2k-2p-1)(2n-2p+1)+2p(2k-2p+1)\Big)
  \label{eq:relreca} \nonumber \\ \\
  &&\!\!\!\!\!\!\!\!\! 
  y_k^{(n)} = -\frac 1 k \sum_{p=0}^{k-1} y_p^{(n)} 
  \frac{(2k+2)!}{(2p+2)!(2k-2p+2)!} \Big((k-p)(2n-2p+1)+p(2k-2p+2)\Big)
  \label{eq:relrecab} \nonumber \\
\end{eqnarray}
with
\begin{equation}
  y_k^{(n)} = \frac{(2k+2)!(2n-2k)!}{2^k k! (2n+1)!} ~ x_k^{(n)}
  \label{eq:xy}
\end{equation}
and the ``initial condition'' $y_0^{(n)} = \frac{2}{2n+1}$.
\\
Consider the second 
recursion relation (\ref{eq:relrecab}). 
It is solved using the generating function:
\begin{equation}
  F_n(u) = \sum_{p=0}^\infty ~ \frac{4^p y_p^{(n)}}{(2p+2)!} ~ u^{2p}
  \label{eq:fctgene}
\end{equation}
Indeed, multiplying equation (\ref{eq:relrecab}) by $(2u)^{2k}$ and
summing over $k$, one transforms the recursion relation into the
following differential equation for the function $F_n$:
\begin{equation}
  F'_n(u) = (2n+1) \left( 
  \frac 1u - \frac{\sinh(2u)}{\cosh(2u)-1} \right) F_n(u) 
  \label{eq:equadiff}
\end{equation}
whose solution is
\begin{equation}
  F_n(u) = \frac{1}{2n+1} \left( \frac{u}{\sinh u} \right)^{2n+1}
  \label{eq:solequadiff}
\end{equation}
so that
\begin{equation}
  x^{(n)}_k = 2^{-k} k! {2n \choose 2k} 
  \left. 
    \frac{d^{2k}}{d u^{2k}} \left(\frac u{\sinh u}\right)^{2n+1}
  \right|_{u=0}
  \label{eq:solx}
\end{equation}

The same procedure applies for the first recursion relation
(\ref{eq:relreca}). One finds exactly the same solution for the
coefficients $y_k^{(n)}$. Of course, this is not surprising since it
follows from the existence of $v_0$ proved in the previous section,
although the two recursion relations look different.

For illustration, we give hereafter the Scasimirs of $\cU(osp(1|2n))$,
for $n=1,\dots,5$:
\begin{eqnarray}
  \Sc^{(1)} & = & A^{(1)}_2 - \frac{1}{2} A^{(1)}_0 \nonumber\\
  & = & [1\bar1]-\frac{1}{2} = \sigma_1\sigma_{\bar 1}-\sigma_{\bar
    1}\sigma_1-\frac{1}{2} \nonumber\\
  \Sc^{(2)} & = & A^{(2)}_4 - 5 A^{(2)}_2 + \frac{9}{2} A^{(2)}_0
  \nonumber\\
  & = & [1\bar 1 2\bar 2] - 5\left([1 \bar 1]+[2\bar 2]\right) +
  \frac{9}{2}  \nonumber\\
  \Sc^{(3)} & = & A^{(3)}_6 - \frac{35}{2} A^{(3)}_4 + \frac{259}{2}
  A^{(3)}_2  - \frac{675}{4} A^{(3)}_0 \nonumber\\
  & = & [1\bar 1 2\bar 2 3\bar 3]  - \frac{35}{2} \left([1 \bar
    1 2\bar 2]+[1\bar 1 3\bar 3]+[2\bar 2 3\bar 3] \right)
    + \frac{259}{2}\left([1 \bar 1]+[2\bar 2]+[3\bar 3] \right) -
    \frac{675}{4}  \nonumber\\
  \Sc^{(4)} & = & A^{(4)}_8 - 42 A^{(4)}_6 + 987 A^{(4)}_4 - 9687
                 A^{(4)}_2 + \frac{33075}{2} A^{(4)}_0 \nonumber\\
  \Sc^{(5)} & = & A^{(5)}_{10} - \frac{165}{2} A^{(5)}_8 + 4389 A^{(5)}_6 
  - \frac{259215}{2} A^{(5)}_4 + \frac{3171663}{2} A^{(5)}_2 -
  \frac{13395375}{4}  A^{(5)}_0 \nonumber\\ 
  \label{eq:scasimirs}
\end{eqnarray}

\subsection{On a funny sum rule}

The superalgebra $osp(1|2n)$ can be realized in terms of bosonic
oscillators. This realization is implemented by imposing the
supplementary relations
\begin{equation}
  {}[\sigma_a , \sigma_b] = \frac12 \; g_{ab} \;.
  \label{eq:realization}
\end{equation}
It is known that the centre is trivial in that case. Hence, the set of
operators anticommuting with all the $\sigma_a$ is $\{0\}$. In
particular, the Scasimir operator vanishes. Moreover, it is easy to
compute $A^{(n)}_{2k}$ and one finds  
$A^{(n)}_{2k} =  2^{-k} \; {\displaystyle {n! \over (n-k)!}}\;$. 
This implies the following sum rule among the coefficients
$x^{(n)}_k$: 
\begin{equation}
  \sum_{k=0}^n  ~ 2^{k-n} ~ {n! \over k!} ~ x^{(n)}_k = 0
  \label{eq:sum_rule}
\end{equation}
which amounts, using (\ref{eq:solx}), to 
\begin{equation}
  \sum_{k=0}^n  ~ {2n \choose 2k} 
  \left. 
    \frac{d^{2k}}{d u^{2k}} \left(\frac u{\sinh u}\right)^{2n+1}
  \right|_{u=0}
  = 
  \left. 
    e^{-u} \frac{d^{2n}}{du^{2n}} e^u 
    \left(\frac u{\sinh u}
    \right)^{2n+1}
  \right|_{u=0} = 0 \;.
  \label{eq:sum_rule2}
\end{equation}
This formula can be proved directly using Cauchy's residue
theorem. Indeed, choosing a small contour around the origin, the last
expression becomes
\begin{equation}
  \oint \frac{du}{u^{2n+1}} e^u  
    \left(\frac u{\sinh u}
    \right)^{2n+1} = 0  
    \qquad\mbox{for } n\geq 1 \;.
  \label{eq:contour}
\end{equation}
Using the variable $t=\sinh u$, this follows from the fact that the
function $\displaystyle{\frac{e^u}{\cosh u}} -1$ is an odd function of
$t$. 

\subsection{On the square of the Scasimir}

The square of the Scasimir is obviously a Casimir element. We shall
now characterize this Casimir operator. It is known that
Harish-Chandra's theorem holds for $osp(1|2n)$ \cite{KacSuperRep},
i.e. there exists an isomorphism between the space of Weyl invariant
polynomials in the Cartan generators and the centre $\cZ$ of
$\cU(osp(1|2n))$.  Using a Poincar{\'e}--Birkhoff--Witt basis, any
element $x$ can be written 
$x=x_0+x_1$ with $x_0\in \cU(\cH)$ and
$x_1\in\cN^-\,\cU(osp(1|2n)) + \cU(osp(1|2n))\,\cN^+$ 
where $\cN^- \oplus \cH \oplus
\cN^+$ is the standard Borel decomposition of $osp(1|2n)$.
Let $h$ be the projection $x\mapsto x_0$ within this direct sum
and $\bar h$ its restriction to the centre $\cZ$. 
The Harish-Chandra 
isomorphism is expressed as $\gamma^{-1} \circ \bar h$ where
$\gamma$ is the automorphism of the ring
of the polynomials in the Cartan generators defined by $\gamma:
H_\alpha\mapsto H_\alpha - (n-\alpha+\frac12)$.  

Using the explicit expression of the Scasimir, we deduce that
$\gamma^{-1}\circ \bar h(\Sc^2)$ is a polynomial in the Cartan
generators of degree $2n$, the power of each being at most equal to two
since each $\sigma_a$ appears only once in $\Sc$. Moreover, this
polynomial being symmetric, one has
\begin{equation}
  \gamma^{-1}\circ \bar h(\Sc^2) = \prod_{\alpha=1}^n H_\alpha^2 \;.
  \label{eq:hc}
\end{equation}
The Casimir operator $\Sc^2$ is then equal to 
$\bar h^{-1}\circ \gamma \left(\prod_{\alpha=1}^n H_\alpha^2\right)$. 
\\
As a buy-product, one also has 
\begin{equation}
  \gamma^{-1}\circ h(\Sc) = \prod_{\alpha=1}^n H_\alpha \;.
  \label{eq:hc^2}
\end{equation}
This last relation appears to have a direct generalization to the
quantum case, whereas the generalization of (\ref{eq:v0suite}) is much
less obvious. 

\medskip

{\bf Note added:} We thank Prof. I. Musson  for sending us his
preprints on superalgebras. One of them, ``On the center of the
enveloping algebra of a classical simple Lie superalgebra'' (to appear
in J. Algebra) contains, among other results, the proof of existence
of the Scasimir for $osp(1|2n)$ and a formula equivalent to
(\ref{eq:hc^2}).

\end{document}